\documentclass[thmsa,12pt]{article}%
\usepackage{amsmath}
\usepackage{amssymb}
\usepackage{sw20lart}
\usepackage{amsfonts}
\usepackage{graphicx}%
\begin{document}

\title{Covariant Derivatives of Multivector and Extensor Fields }
\author{{\footnotesize A. M. Moya}$^{2}${\footnotesize , V. V. Fern\'{a}ndez}$^{1}%
${\footnotesize and W. A. Rodrigues Jr}.$^{1}${\footnotesize \ }\\$^{1}\hspace{-0.1cm}${\footnotesize Institute of Mathematics, Statistics and
Scientific Computation}\\{\footnotesize \ IMECC-UNICAMP CP 6065}\\{\footnotesize \ 13083-859 Campinas, SP, Brazil }\\{\footnotesize e-mail:walrod@ime.unicamp.br}\\{\footnotesize \ }$^{2}${\footnotesize Department of Mathematics, University
of Antofagasta, Antofagasta, Chile} \\{\footnotesize e-mail: mmoya@uantof.cl}}
\maketitle

\begin{abstract}
We give in this paper which is the fifth in a series of eight a theory of
covariant derivatives of multivector and extensor fields based on the
geometric calculus of an arbitrary smooth manifold $M$, and the notion of a
connection extensor field defining a parallelism structure on $M.$ Also we
give a novel and intrinsic presentation (i.e., one that does not depend on a
chosen orthonormal moving frame) of the torsion and curvature fields of
Cartan's theory. Two kinds of Cartan's connection fields are identified, and
both appear in the intrinsic Cartan's structure equations satisfied by the
Cartan's torsion and curvature extensor fields.

\end{abstract}
\tableofcontents

\section{Introduction}

This is the fifth paper in a series of eight. It is dedicated to a theory of
covariant derivatives of multivector and extensor fields. We introduce in
Section 2 the notion of a \emph{parallelism structure} on a smooth manifold
$M$ given by a \emph{connection extensor field} $\gamma$.

In Sections 3 and 4 we introduce the concepts of $a$-directional covariant
derivatives of multivector and extensor fields, respectively, and prove the
main properties satisfied by these objects. In Section 5 we give a thoughtful
study of the so-called \emph{symmetric parallelism structures}, where among
others we present and prove a \emph{Bianchi-like identity} and give an
intrinsic Cartan theory of the torsion and curvature fields. \emph{Cartan's
connections of the first and second kind }are identified, and both appear on
Cartan's structure equations. We emphasize the novelty of our approach to
Cartan theory, namely that it does not depend on a chosen orthonormal moving
frame, hence, the name\textit{ intrinsic} used above. In an Appendix several
examples are worked in detail in order to show how our theory relates to known
theories dealing with the same subject.

\section{Parallelism Structure}

Let $M$ be a $n$-dimensional smooth manifold. Then, for any point $o\in M$
there exists a local coordinate system $(U_{o},\phi_{o})$ such that $o\in
U_{o}$ and $\phi_{o}(o)=(0,\ldots,0)\in\mathbb{R}^{n}.$ $U_{o}$ is an open
subset of $M,$ and $\phi_{o}$ is an homeomorphism from $U_{o}$ onto the open
subset $\phi_{o}(U_{o})\subseteq\mathbb{R}^{n}.$

As in \cite{1}, let $\mathcal{U}_{o}$ be the \emph{canonical space} for
$(U_{o},\phi_{o})$, and $U$ be an open subset of $U_{o}.$ We denote the ring
(with identity) of smooth scalar fields on $U,$ the module of smooth vector
fields on $U$ and the module of smooth multivector fields on $U$ respectively
by $\mathcal{S}(U),$ $\mathcal{V}(U)$ and $\mathcal{M}(U).$ The set of smooth
$k$-vector fields on $U$ is denoted by $\mathcal{M}^{k}(U).$ The module of
smooth $k$-extensor fields on $U$ is denoted by $k$-$ext(\mathcal{M}%
_{1}^{\diamond}(U),\ldots,\mathcal{M}_{k}^{\diamond}(U);\mathcal{M}^{\diamond
}(U))$, see \cite{1}.

Any smooth vector elementary $2$-extensor field on $U$ is said to be a
\emph{connection field} on $U.$ A general connection field will be denoted by
$\gamma,$ i.e., $\gamma:U\rightarrow2$-$ext^{1}(\mathcal{U}_{o})$. The
smoothness of such $\gamma$ means that for all $a,b\in\mathcal{V}%
(\mathcal{U}_{o})$ the vector field defined by $U\ni p\mapsto\gamma
_{(p)}(a(p),b(p))\in\mathcal{U}_{o}$ is itself smooth.

The open set $U$ equipped with such a connection field $\gamma,$ namely
$(U,\gamma),$ will be said to be a \emph{parallelism structure} on $U.$

Let us take $a\in\mathcal{U}_{o}.$ The smooth $(1,1)$-extensor field on $U,$
namely $\gamma_{a},$ defined as $U\ni p\mapsto\left.  \gamma_{a}\right|
_{(p)}\in ext_{1}^{1}(\mathcal{U}_{o})$ such that for all $b\in\mathcal{U}%
_{o}$
\begin{equation}
\left.  \gamma_{a}\right|  _{(p)}(b)=\gamma_{(p)}(a,b), \label{PS.1}%
\end{equation}
will be called the $a$-\emph{directional connection field on }$U$, of course,
associated to $(U,\gamma).$

We emphasize that the $(1,1)$-extensor character and the smoothness of
$\gamma_{a}$ are immediate consequences of the vector elementary $2$-extensor
character and the smoothness of $\gamma.$

The smooth $(1,2)$-extensor field on $U,$ namely $\Omega,$ defined as $U\ni
p\mapsto\Omega_{(p)}\in ext_{1}^{2}(\mathcal{U}_{o})$ such that for all
$a\in\mathcal{U}_{o}$
\begin{equation}
\Omega_{(p)}(a)=\dfrac{1}{2}biv[\left.  \gamma_{a}\right|  _{(p)}],
\label{PS.2}%
\end{equation}
will be called the \emph{gauge connection field on }$U.$

From the definition of $biv[t]$ (see \cite{2}), taking any pair of reciprocal
frame fields on $U,$ say $(\{e_{\mu}\},\{e^{\mu}\}),$ and using Eq.(\ref{PS.1}%
), we can write Eq.(\ref{PS.2}) as
\begin{equation}
\Omega_{(p)}(a)=\frac{1}{2}\gamma_{(p)}(a,e^{\mu}(p))\wedge e_{\mu}%
(p)=\frac{1}{2}\gamma_{(p)}(a,e_{\mu}(p))\wedge e^{\mu}(p). \label{PS.2a}%
\end{equation}
So, we see that the $(1,2)$-extensor character and the smoothness of $\Omega$
are easily deduced from the vector elementary $2$-extensor character and the
smoothness of $\gamma.$

Let us take any pair of reciprocal frame fields on $U,$ say $(\{e_{\mu
}\},\{e^{\mu}\})$, i.e., $e_{\mu}\cdot e^{\nu}=\delta_{\mu}^{\nu}.$ Let
$\Gamma_{a}$ be the \emph{generalized (extensor field)} of $\gamma_{a}$ (see
\cite{2}), i.e., $\Gamma_{a}$ defined as $U\ni p\mapsto\left.  \Gamma
_{a}\right\vert _{(p)}\in ext(\mathcal{U}_{o}),$ is a smooth extensor field on
$U$ such that for all $X\in\bigwedge\mathcal{U}_{o}$
\begin{equation}
\left.  \Gamma_{a}\right\vert _{(p)}(X)=\left.  \gamma_{a}\right\vert
_{(p)}(e^{\mu}(p))\wedge(e_{\mu}(p)\lrcorner X)=\left.  \gamma_{a}\right\vert
_{(p)}(e_{\mu}(p))\wedge(e^{\mu}(p)\lrcorner X). \label{PS.3}%
\end{equation}

It is easily seen that (for a given $X$) the multivector appearing on the
right side of Eq.(\ref{PS.3}) does not depend on the choice of the reciprocal
frame fields. The extensor character and the smoothness of $\Gamma_{a}$
follows from the $(1,1)$-extensor character and the smoothness of $\gamma
_{a}.$

We will usually omit the letter $p$ in writing the definitions given by
Eq.(\ref{PS.1}), Eq.(\ref{PS.2}) and Eq.(\ref{PS.3}), and other equations
using extensor fields. No confusion should arise with this standard practice.
Eq.(\ref{PS.3}) might be also written in the succinct form $\Gamma
_{a}(X)=\gamma_{a}(\partial_{b})\wedge(b\lrcorner X).$

We end this section by presenting some of the basic properties satisfied by
$\Gamma_{a}$.\vspace{0.1in}

\textbf{i.} $\Gamma_{a}$ is grade-preserving, i.e.,
\begin{equation}
\text{if }X\in\mathcal{M}^{k}(U),\text{ then }\Gamma_{a}(X)\in\mathcal{M}%
^{k}(U). \label{PS.4}%
\end{equation}

\textbf{ii}. For any $X\in\mathcal{M}(U)$,
\begin{align}
\Gamma_{a}(\widehat{X})  &  =\widehat{\Gamma_{a}(X)},\label{PS.5a}\\
\Gamma_{a}(\widetilde{X})  &  =\widetilde{\Gamma_{a}(X)},\label{PS.5b}\\
\Gamma_{a}(\overline{X})  &  =\overline{\Gamma_{a}(X)}. \label{PS.5c}%
\end{align}

\textbf{iii.} For any $f\in\mathcal{S}(U),$ $b\in\mathcal{V}(U)$ and
$X,Y\in\mathcal{M}(U),$
\begin{align}
\Gamma_{a}(f)  &  =0,\label{PS.6a}\\
\Gamma_{a}(b)  &  =\gamma_{a}(b),\label{PS.6b}\\
\Gamma_{a}(X\wedge Y)  &  =\Gamma_{a}(X)\wedge Y+X\wedge\Gamma_{a}(Y).
\label{PS.6c}%
\end{align}

\textbf{iv.} The adjoint of $\Gamma_{a},$ namely $\Gamma_{a}^{\dagger},$ is
the generalized of the adjoint of $\gamma_{a}$, namely $\gamma_{a}^{\dagger},$
i.e.,
\begin{equation}
\Gamma_{a}^{\dagger}(X)=\gamma_{a}^{\dagger}(\partial_{b})\wedge(b\lrcorner
X). \label{PS.7}%
\end{equation}

\textbf{v}. The symmetric (skew-symmetric) part of $\Gamma_{a},$ namely
$\Gamma_{a\pm}=\dfrac{1}{2}(\Gamma_{a}\pm\Gamma_{a}^{\dagger}),$ is the
generalized of the symmetric (skew-symmetric) part of $\gamma_{a},$ namely
$\gamma_{a\pm}=\dfrac{1}{2}(\gamma_{a}\pm\gamma_{a}^{\dagger}),$ i.e.,
\begin{equation}
\Gamma_{a\pm}(X)=\gamma_{a\pm}(\partial_{b})\wedge(b\lrcorner X). \label{PS.8}%
\end{equation}

\textbf{vi}. $\Gamma_{a-}$ can be factorized by a remarkable formula which
only involves $\Omega.$ It is
\begin{equation}
\Gamma_{a-}(X)=\Omega(a)\times X. \label{PS.9}%
\end{equation}

\textbf{vii}. For any $X,Y\in\mathcal{M}(U)$ it holds
\begin{equation}
\Gamma_{a-}(X\ast Y)=\Gamma_{a-}(X)\ast Y+X\ast\Gamma_{a-}(Y), \label{PS.10}%
\end{equation}
where $\ast$ means any suitable product of smooth multivector fields, either
$(\wedge),$ $(\cdot),$ $(\lrcorner,\llcorner)$ or $(b$-\emph{Clifford
product}$),$ see \cite{1,3}.

\section{$a$-Directional Covariant Derivatives of Multivector Fields}

Given a parallelism structure $(U,\gamma),$ let us take $a\in\mathcal{U}_{o}$.
Then, associated to $(U,\gamma)$ we can introduce two $a$-\emph{directional
covariant derivative operators} ($a$-\emph{DCDO's}), namely $\nabla_{a}^{+}$
and $\nabla_{a}^{-},$ which act on the module of smooth multivector fields on
$U.$

They are defined by $\nabla_{a}^{\pm}:\mathcal{M}(U)\rightarrow\mathcal{M}(U)$
such that
\begin{align}
\nabla_{a}^{+}X(p)  &  =a\cdot\partial_{o}X(p)+\left.  \Gamma_{a}\right|
_{(p)}(X(p))\label{CDM.1a}\\
\nabla_{a}^{-}X(p)  &  =a\cdot\partial_{o}X(p)-\left.  \Gamma_{a}^{\dagger
}\right|  _{(p)}(X(p))\text{ for each }p\in U, \label{CDM.1b}%
\end{align}
where $a\cdot\partial_{o}$ is the canonical $a$-\emph{DODO} as was defined in
\cite{1}.

We emphasize that each of $\nabla_{a}^{+}$ and $\nabla_{a}^{-}$ satisfies
\emph{indeed }the fundamental properties which a well-defined covariant
derivative is expected to have. This is trivial to verify whenever we take
into account the well-known properties of $a\cdot\partial_{o}$ (see \cite{1}),
and the properties of $\Gamma_{a}$ given by Eq.(\ref{PS.4}), Eqs.(\ref{PS.6a}%
), (\ref{PS.6b}) and (\ref{PS.6c}), and Eq.(\ref{PS.7}).

As usual we will write Eq.(\ref{CDM.1a}) and Eq.(\ref{CDM.1b}) by omitting $p
$ when no confusion arises.

The smooth multivector fields on $U,$ namely $\nabla_{a}^{+}X$ and $\nabla
_{a}^{-}X,$ will be respectively called the \emph{plus }and the \emph{minus}
$a$-\emph{directional covariant derivatives} of $X.$

We summarize some of the most important properties for the pair of
$a$-\emph{DCDO's} $\nabla_{a}^{+}$ and $\nabla_{a}^{-}.\vspace{0.1in}$

\textbf{i.} $\nabla_{a}^{+}$ and $\nabla_{a}^{-}$ are grade-preserving
operators on $\mathcal{M}(U),$ i.e.,
\begin{equation}
\text{if }X\in\mathcal{M}^{k}(U),\text{ then }\nabla_{a}^{\pm}X\in
\mathcal{M}^{k}(U). \label{CDM.2}%
\end{equation}

\textbf{ii. }For all $X\in\mathcal{M}(U),$ and for any $\alpha,\alpha^{\prime
}\in\mathbb{R}$ and $a,a^{\prime}\in\mathcal{U}_{o}$ we have
\begin{equation}
\nabla_{\alpha a+\alpha^{\prime}a^{\prime}}^{\pm}X=\alpha\nabla_{a}^{\pm
}X+\alpha^{\prime}\nabla_{a^{\prime}}^{\pm}X. \label{CDM.3}%
\end{equation}

\textbf{iii. }For all $f\in\mathcal{S}(U)$ and $X,Y\in\mathcal{M}(U)$ we have
\begin{align}
\nabla_{a}^{\pm}f  &  =a\cdot\partial_{o}f,\label{CDM.4a}\\
\nabla_{a}^{\pm}(X+Y)  &  =\nabla_{a}^{\pm}X+\nabla_{a}^{\pm}Y, \label{CDM.4b}%
\\
\nabla_{a}^{\pm}(fX)  &  =(a\cdot\partial_{o}f)X+f(\nabla_{a}^{\pm}X).
\label{CDM.4c}%
\end{align}

\textbf{iv. }For all $X,Y\in\mathcal{M}(U)$ we have
\begin{equation}
\nabla_{a}^{\pm}(X\wedge Y)=(\nabla_{a}^{\pm}X)\wedge Y+X\wedge(\nabla
_{a}^{\pm}Y). \label{CDM.5}%
\end{equation}

\textbf{v. }For all $X,Y\in\mathcal{M}(U)$ we have
\begin{equation}
(\nabla_{a}^{+}X)\cdot Y+X\cdot(\nabla_{a}^{-}Y)=a\cdot\partial_{o}(X\cdot Y).
\label{CDM.6}%
\end{equation}

It should be noticed that $(\nabla_{a}^{+},\nabla_{a}^{-})$ as defined by
Eq.(\ref{CDM.1a}) and Eq.(\ref{CDM.1b}) is the unique pair of $a$%
-\emph{DCDO's} associated to $(U,\gamma)$ which satisfies the remarkable
property given by Eq.(\ref{CDM.6}).

We emphasize here that the $a$-\emph{DODO }$a\cdot\partial_{o}$ acting on
$\mathcal{M}(U),$ is also a well-defined $a$-\emph{DCDO.} In this particular
case, the connection field $\gamma$ is \emph{identically zero} and the plus
and minus $a$-\emph{DCDO's} are equal to each other, and both of them coincide
with $a\cdot\partial_{o}.$

We introduce yet another well-defined $a$-\emph{DCDO} which acts also on the
module of smooth multivector fields on $U,$ namely $\nabla_{a}^{0}.$

It is defined by
\begin{equation}
\nabla_{a}^{0}X=\frac{1}{2}(\nabla_{a}^{+}X+\nabla_{a}^{-}X). \label{CDM.7}%
\end{equation}
But, by using Eqs.(\ref{CDM.1a}) and (\ref{CDM.1b}), and Eq.(\ref{PS.9}), we
might write else
\begin{equation}
\nabla_{a}^{0}X=a\cdot\partial_{o}X+\Omega(a)\times X. \label{CDM.8}%
\end{equation}

The $a$-\emph{DCDO} $\nabla_{a}^{0}$ satisfies the same properties which hold
for each one of the $a$-\emph{DCDO's} $\nabla_{a}^{+}$ and $\nabla_{a}^{-}.$
But, it has also an additional remarkable property
\begin{equation}
(\nabla_{a}^{0}X)\cdot Y+X\cdot(\nabla_{a}^{0}Y)=a\cdot\partial_{o}(X\cdot Y).
\label{CDM.9}%
\end{equation}

Moreover, it satisfies a Leibnitz-like rule for any suitable product of smooth
multivector fields, i.e.,
\begin{equation}
\nabla_{a}^{0}(X*Y)=(\nabla_{a}^{0}X)*Y+X*(\nabla_{a}^{0}Y). \label{CDM.10}%
\end{equation}

\subsection{Connection Operators}

Associated to any parallelism structure $(U,\gamma)$ we can introduce two
remarkable operators which map $2$-uples of smooth vector fields to smooth
vector fields.

They are defined by $\Gamma^{\pm}:\mathcal{V}(U)\times\mathcal{V}%
(U)\rightarrow\mathcal{V}(U)$ such that
\begin{equation}
\Gamma^{\pm}(a,b)=\nabla_{a}^{\pm}b, \label{CO.1}%
\end{equation}
and will be called the \emph{connection operators} of $(U,\gamma).$

We summarize the basic properties of them.\vspace{0.1in}

\textbf{i.} For all $f\in\mathcal{S}(U),$ and $a,a^{\prime},b,b^{\prime}%
\in\mathcal{V}(U),$ we have
\begin{align}
\Gamma^{\pm}(a+a^{\prime},b)  &  =\Gamma^{\pm}(a,b)+\Gamma^{\pm}(a^{\prime
},b),\label{CO.2a}\\
\Gamma^{\pm}(a,b+b^{\prime})  &  =\Gamma^{\pm}(a,b)+\Gamma^{\pm}(a,b^{\prime
}),\label{CO.2b}\\
\Gamma^{\pm}(fa,b)  &  =f\Gamma^{\pm}(a,b),\label{CO.2c}\\
\Gamma^{\pm}(a,fb)  &  =(a\cdot\partial_{o}f)b+f\Gamma^{\pm}(a,b).
\label{CO.2d}%
\end{align}
As we can observe both connection operators satisfy the linearity property
only with respect to the first smooth vector field variable. Thus, connection
operators are not extensor fields.

\textbf{ii.} For all $a,b,c\in\mathcal{V}(U),$ we have
\begin{equation}
\Gamma^{+}(a,b)\cdot c+b\cdot\Gamma^{-}(a,c)=a\cdot\partial_{o}(b\cdot c).
\label{CO.3}%
\end{equation}
It is an immediate consequence of Eq.(\ref{CDM.6}).

\subsection{Deformation of Covariant Derivatives}

Let $(\nabla_{a}^{+},\nabla_{a}^{-})$ be any pair of $a$-\emph{DCDO's} \ and
$\lambda$ a non-singular smooth $(1,1)$-extensor field on $U$ $.$ We define
the deformation of these covariant derivatives as the pair $(_{\lambda}%
\nabla_{a}^{+},_{\lambda}\nabla_{a}^{-})$ by
\begin{align}
_{\lambda}\nabla_{a}^{+}X  &  =\underline{\lambda}(\nabla_{a}^{+}%
\underline{\lambda}^{-1}(X)),\\
_{\lambda}\nabla_{a}^{-}X  &  =\underline{\lambda}^{\ast}(\nabla_{a}%
^{-}\underline{\lambda}^{\dagger}(X)),
\end{align}
where $\underline{\lambda}$ is the \emph{extended}\footnote{Recall that
$\lambda^{*}=(\lambda^{-1})^{\dagger}=(\lambda^{\dagger})^{-1},$ and
$\underline{\lambda}^{-1}=(\underline{\lambda})^{-1}=\underline{(\lambda
^{-1})}$ and $\underline{\lambda}^{\dagger}=(\underline{\lambda})^{\dagger
}=\underline{(\lambda^{\dagger})},$ see \cite{2}$.$} of $\lambda$, is a
well-defined pair of $a$-\emph{DCDO's,} since it satisfies as it is trivial to
verify the fundamental properties given by Eqs.(\ref{CDM.4a}), (\ref{CDM.4b})
and (\ref{CDM.4c}), Eq.(\ref{CDM.5}) and Eq.(\ref{CDM.6}). For instance,%
\begin{equation}
_{\lambda}\nabla_{a}^{+}f=\underline{\lambda}(\nabla_{a}^{+}\underline
{\lambda}^{-1}(f))=a\cdot\partial_{0}f.
\end{equation}
\ \ 

We verify now that the definitions given by Eq.(\ref{CDM.11a}) and
Eq.(\ref{CDM.11b}) also satisfy a property analogous to the one given by
Eq.(\ref{CDM.6}). Indeed, we have
\begin{align}
(_{\lambda}\nabla_{a}^{+}X)\cdot Y+X\cdot(_{\lambda}\nabla_{a}^{-}Y)  &
=(\nabla_{a}^{+}\underline{\lambda}^{-1}(X))\cdot\underline{\lambda}^{\dagger
}(Y)+\underline{\lambda}^{-1}(X)\cdot(\nabla_{a}^{-}\underline{\lambda
}^{\dagger}(Y))\nonumber\\
&  =a\cdot\partial_{o}(\underline{\lambda}^{-1}(X)\cdot\underline{\lambda
}^{\dagger}(Y)),\nonumber\\
&  =a\cdot\partial_{o}(X\cdot Y).
\end{align}

\section{$a$-Directional Covariant Derivatives of Extensor Fields}

The three $a$-\emph{DCDO's} $\nabla_{a}^{+},$ $\nabla_{a}^{-}$ and $\nabla
_{a}^{0}$ which act on $\mathcal{M}(U)$ can be extended in order to act on the
module of smooth $k$-extensor fields on $U$. For any $t\in k$-$ext(\mathcal{M}%
_{1}^{\diamond}(U),\ldots,\mathcal{M}_{k}^{\diamond}(U);\mathcal{M}^{\diamond
}(U))$, we can define exactly $3^{k+1}$ covariant derivatives, namely
$\nabla_{a}^{\sigma_{1}\ldots\sigma_{k}\sigma}t\in k $-$ext(\mathcal{M}%
_{1}^{\diamond}(U),\ldots,\mathcal{M}_{k}^{\diamond}(U);\mathcal{M}^{\diamond
}(U)),$ where each of $\sigma_{1},\ldots,\sigma_{k},\sigma$ is being used to
mean either $(+),$ $(-)$ or $(0).$ They are given by the following definition.

For all $X_{1}\in\mathcal{M}_{1}^{\diamond}(U),\ldots,X_{k}\in\mathcal{M}%
_{k}^{\diamond}(U),$ $X\in\mathcal{M}^{\diamond}(U)$
\begin{align}
&  (\nabla_{a}^{\sigma_{1}\ldots\sigma_{k}\sigma}t)_{(p)}(X_{1}(p),\ldots
,X_{k}(p))\cdot X(p)\nonumber\\
&  =a\cdot\partial_{o}(t_{(p)}(\ldots)\cdot X(p))-t_{(p)}(\nabla_{a}%
^{\sigma_{1}}X_{1}(p),\ldots)\cdot X(p)\nonumber\\
&  -\cdots-t_{(p)}(\ldots,\nabla_{a}^{\sigma_{k}}X_{k}(p))\cdot X(p)-t_{(p)}%
(\ldots)\cdot\nabla_{a}^{\sigma}X(p), \label{CDE.1}%
\end{align}
for each $p\in U$.

As usual when no confusion arises we will write Eq.(\ref{CDE.1}) by omitting
$p.$

We call the reader's attention that each of $\nabla_{a}^{\sigma_{1}%
\ldots\sigma_{k}\sigma}t$ as defined by Eq.(\ref{CDE.1}) is in fact a smooth
$k$-extensor field. Its $k$-extensor character and smoothness can be easily
deduced from the respective properties of $t.$ We note also that in the first
term on the right side of Eq.(\ref{CDE.1}), $a\cdot\partial_{o}$ refers to the
canonical $a$-\emph{DODO} as was defined in \cite{1}.

We notice that any smooth $(1,1)$-extensor field on $U,$ say $t,$ has just
$3^{1+1}=9$ covariant derivatives. For instance, four important covariant
derivatives of such $t$ are given by
\begin{align}
(\nabla_{a}^{++}t)(X_{1})\cdot X  &  =a\cdot\partial_{o}(t(X_{1})\cdot
X)\nonumber\\
&  -t(\nabla_{a}^{+}X_{1})\cdot X-t(X_{1})\cdot\nabla_{a}^{+}X, \label{CDE.1a}%
\\
(\nabla_{a}^{+-}t)(X_{1})\cdot X  &  =a\cdot\partial_{o}(t(X_{1})\cdot
X)\nonumber\\
&  -t(\nabla_{a}^{+}X_{1})\cdot X-t(X_{1})\cdot\nabla_{a}^{-}X, \label{CDE.1b}%
\\
(\nabla_{a}^{--}t)(X_{1})\cdot X  &  =a\cdot\partial_{o}(t(X_{1})\cdot
X)\nonumber\\
&  -t(\nabla_{a}^{-}X_{1})\cdot X-t(X_{1})\cdot\nabla_{a}^{-}X, \label{CDE.1c}%
\\
(\nabla_{a}^{-+}t)(X_{1})\cdot X  &  =a\cdot\partial_{o}(t(X_{1})\cdot
X)\nonumber\\
&  -t(\nabla_{a}^{-}X_{1})\cdot X-t(X_{1})\cdot\nabla_{a}^{+}X, \label{CDE.1d}%
\end{align}
where $X_{1}\in\mathcal{M}_{1}^{\diamond}(U)$ and $X\in\mathcal{M}^{\diamond
}(U),$

We present now some of the basic properties satisfied by these $a$-directional
covariant derivatives of smooth $k$-extensor fields.\vspace{0.1in}

\textbf{i. }For all $f\in\emph{S}(U),$ and $t,u\in k$-$ext(\mathcal{M}%
_{1}^{\diamond}(U),\ldots,\mathcal{M}_{k}^{\diamond}(U);\mathcal{M}^{\diamond
}(U)),$ it holds
\begin{align}
\nabla_{a}^{\sigma_{1}\ldots\sigma_{k}\sigma}(t+u)  &  =\nabla_{a}^{\sigma
_{1}\ldots\sigma_{k}\sigma}t+\nabla_{a}^{\sigma_{1}\ldots\sigma_{k}\sigma
}u,\label{CDE.2a}\\
\nabla_{a}^{\sigma_{1}\ldots\sigma_{k}\sigma}(ft)  &  =(a\cdot\partial
_{o}f)t+f(\nabla_{a}^{\sigma_{1}\ldots\sigma_{k}\sigma}t). \label{CDE.2b}%
\end{align}

\textbf{ii.} For all $t\in1$-$ext(\mathcal{M}_{1}^{\diamond}(U);\mathcal{M}%
^{\diamond}(U)),$ it holds a noticeable property
\begin{equation}
(\nabla_{a}^{\sigma_{1}\sigma}t)^{\dagger}=\nabla_{a}^{\sigma\sigma_{1}%
}t^{\dagger}. \label{CDE.3}%
\end{equation}
Note the inversion between $\sigma_{1}$ and $\sigma$ into the $a$-\emph{DCDO's
}above. As we can see, the three $a$-\emph{DCDO's }$\nabla_{a}^{++},$
$\nabla_{a}^{--}$ and $\nabla_{a}^{00}$ commute indeed with the adjoint
operator $\left.  {}\right.  ^{\dagger}.$

\textbf{Proof}

Let us take $X_{1}\in\mathcal{M}_{1}^{\diamond}(U)$ and $X\in\mathcal{M}%
^{\diamond}(U).$ By recalling the fundamental property of the adjoint operator
\cite{2}, and in accordance with Eq.(\ref{CDE.1}), we can write
\begin{align*}
(\nabla_{a}^{\sigma_{1}\sigma}t)^{\dagger}(X)\cdot X_{1}  &  =(\nabla
_{a}^{\sigma_{1}\sigma}t)(X_{1})\cdot X\\
&  =a\cdot\partial_{o}(t(X_{1})\cdot X)-t(\nabla_{a}^{\sigma_{1}}X_{1})\cdot
X-t(X_{1})\cdot\nabla_{a}^{\sigma}X\\
&  =a\cdot\partial_{o}(t^{\dagger}(X)\cdot X_{1})-t^{\dagger}(X)\cdot
\nabla_{a}^{\sigma_{1}}X_{1}-t^{\dagger}(\nabla_{a}^{\sigma}X)\cdot X_{1},\\
&  =(\nabla_{a}^{\sigma\sigma_{1}}t^{\dagger})(X)\cdot X_{1}.
\end{align*}
Hence, by non-degeneracy of scalar product, the expected result immediately
follows.$\blacksquare$

\textbf{iii.} For all $t\in1$-$ext(\mathcal{M}_{1}^{\diamond}(U),\mathcal{M}%
^{\diamond}(U)),$ it holds
\begin{align}
(\nabla_{a}^{++}t)(X_{1})  &  =\nabla_{a}^{-}t(X_{1})-t(\nabla_{a}^{+}%
X_{1}),\label{CDE.4a}\\
(\nabla_{a}^{+-}t)(X_{1})  &  =\nabla_{a}^{+}t(X_{1})-t(\nabla_{a}^{+}%
X_{1}),\label{CDE.4b}\\
(\nabla_{a}^{--}t)(X_{1})  &  =\nabla_{a}^{+}t(X_{1})-t(\nabla_{a}^{-}%
X_{1}),\label{CDE.4c}\\
\nabla_{a}^{-+}t)(X_{1})  &  =\nabla_{a}^{-}t(X_{1})-t(\nabla_{a}^{-}X_{1}).
\label{CDE.4d}%
\end{align}

\textbf{Proof}

We will only try Eq.(\ref{CDE.4a}). Let us take $X_{1}\in\mathcal{M}%
_{1}^{\diamond}(U)$ and $X\in\mathcal{M}^{\diamond}(U).$ In accordance with
Eq.(\ref{CDE.1}) and by recalling Eq.(\ref{CDM.6}), we have
\begin{align*}
(\nabla_{a}^{++}t)(X_{1})\cdot X  &  =a\cdot\partial_{o}(t(X_{1})\cdot
X)-t(\nabla_{a}^{+}X_{1})\cdot X-t(X_{1})\cdot\nabla_{a}^{+}X\\
&  =\nabla_{a}^{-}t(X_{1})\cdot X+t(X_{1})\cdot\nabla_{a}^{+}X\\
&  -t(\nabla_{a}^{+}X_{1})\cdot X-t(X_{1})\cdot\nabla_{a}^{+}X\\
&  =(\nabla_{a}^{-}t(X_{1})-t(\nabla_{a}^{+}X_{1}))\cdot X.
\end{align*}
Hence, by non-degeneracy of scalar product, it follows what was to be proved.
$\blacksquare$

\section{Torsion and Curvature Fields}

Let $(U,\gamma)$ be a parallelism structure on $U.$ The \emph{smooth vector
elementary} $2$-\emph{exform field} on $U,$ namely, $\tau$ such that for all
$a,b\in\mathcal{V}(U)$
\begin{equation}
\tau(a,b)=\nabla_{a}^{+}b-\nabla_{b}^{+}a-[a,b], \label{TCF.1a}%
\end{equation}
i.e.,
\begin{equation}
\tau(a,b)=\gamma_{a}(b)-\gamma_{b}(a), \label{TCF.1b}%
\end{equation}
will be called the \emph{torsion field} of $(U,\gamma).$

The \emph{smooth vector elementary} $3$-\emph{extensor field }on $U,$ namely
$\rho,$ such that for all $a,b,c\in\mathcal{V}(U)$
\begin{equation}
\rho(a,b,c)=[\nabla_{a}^{+},\nabla_{b}^{+}]c-\nabla_{[a,b]}^{+}c,
\label{TCF.2a}%
\end{equation}
i.e.,
\begin{equation}
\rho(a,b,c)=(a\cdot\partial_{o}\gamma_{b})(c)-(b\cdot\partial_{o}\gamma
_{a})(c)+[\gamma_{a},\gamma_{b}](c)-\gamma_{[a,b]}(c), \label{TCF.2b}%
\end{equation}
will be called the \emph{curvature field} of $(U,\gamma).$

It should be emphasized that the curvature field $\rho$ is
\emph{skew-symmetric} in first and second variables, i.e.,
\begin{equation}
\rho(a,b,c)=-\rho(b,a,c). \label{TCF.3}%
\end{equation}

\subsection{Symmetric Parallelism Structures}

A parallelism structure $(U,\gamma)$ is said to be \emph{symmetric} if and
only if
\begin{equation}
\gamma_{a}(b)=\gamma_{b}(a). \label{SPS.1}%
\end{equation}
As we can easily prove this condition above is completely equivalent to the
following condition. For all $a,b\in\mathcal{V}(U)$
\begin{equation}
\nabla_{a}^{+}b-\nabla_{b}^{+}a=[a,b]. \label{SPS.2}%
\end{equation}

In accordance with Eq.(\ref{TCF.1a}) and Eq.(\ref{TCF.1b}) we have that a
parallelism structure is symmetric if and only if it is torsionless, i.e.,
\begin{equation}
\tau(a,b)=0. \label{SPS.3}%
\end{equation}

We now present and prove some basic properties of a symmetric parallelism
structure.\vspace{0.1in}

\textbf{i.} The curvature field $\rho$ satisfies the cyclic property
\begin{equation}
\rho(a,b,c)+\rho(b,c,a)+\rho(c,a,b)=0. \label{SPS.4}%
\end{equation}

\textbf{Proof}

By recalling Eq.(\ref{TCF.2a}) we can write
\begin{align}
\rho(a,b,c)  &  =\nabla_{a}^{+}\nabla_{b}^{+}c-\nabla_{b}^{+}\nabla_{a}%
^{+}c-\nabla_{\lbrack a,b]}^{+}c,\label{SPS.4a}\\
\rho(b,c,a)  &  =\nabla_{b}^{+}\nabla_{c}^{+}a-\nabla_{c}^{+}\nabla_{b}%
^{+}a-\nabla_{\lbrack b,c]}^{+}a,\label{SPS.4b}\\
\rho(c,a,b)  &  =\nabla_{c}^{+}\nabla_{a}^{+}b-\nabla_{a}^{+}\nabla_{c}%
^{+}b-\nabla_{\lbrack c,a]}^{+}b. \label{SPS.4c}%
\end{align}

By adding Eqs.(\ref{SPS.4a}), (\ref{SPS.4b}) and (\ref{SPS.4c}), wherever by
taking into account Eq.(\ref{SPS.2}), we get
\begin{align}
&  \rho(a,b,c)+\rho(b,c,a)+\rho(c,a,b)\nonumber\\
&  =\nabla_{a}^{+}(\nabla_{b}^{+}c-\nabla_{c}^{+}b)+\nabla_{b}^{+}(\nabla
_{c}^{+}a-\nabla_{a}^{+}c)+\nabla_{c}^{+}(\nabla_{a}^{+}b-\nabla_{b}%
^{+}a)\nonumber\\
&  -\nabla_{\lbrack a,b]}^{+}c-\nabla_{\lbrack b,c]}^{+}a-\nabla_{\lbrack
c,a]}^{+}b,\nonumber\\
&  =[a,[b,c]]+[b,[c,a]]+[c,[a,b]]. \label{SPS.4d}%
\end{align}
Hence, by recalling the so-called Jacobi's identity for the Lie product of
smooth vector fields \cite{1}, the expected result immediate
follows.$\blacksquare$

\textbf{ii.} The curvature field $\rho$ satisfies the so-called
\emph{Bianchi's identity}, i.e.,
\begin{equation}
(\nabla_{d}^{+++-}\rho)(a,b,c)+(\nabla_{a}^{+++-}\rho)(b,d,c)+(\nabla
_{b}^{+++-}\rho)(d,a,c)=0. \label{SPS.5}%
\end{equation}

\textbf{Proof}

Let us take $a,b,c,d,w\in\mathcal{V}(U).$ In accordance with Eq.(\ref{CDE.1}),
by using Eq.(\ref{CDM.6}), we have
\begin{align*}
&  (\nabla_{d}^{+++-}\rho)(a,b,c)\cdot w\\
&  =d\cdot\partial_{o}(\rho(a,b,c)\cdot w)-\rho(\nabla_{d}^{+}a,b,c)\cdot
w-\rho(a,\nabla_{d}^{+}b,c)\cdot w\\
&  -\rho(a,b,\nabla_{d}^{+}c)\cdot w-\rho(a,b,c)\cdot\nabla_{d}^{-}w,
\end{align*}
i.e.,
\begin{align}
(\nabla_{d}^{+++-}\rho)(a,b,c)  &  =\nabla_{d}^{+}\rho(a,b,c)-\rho(\nabla
_{d}^{+}a,b,c)\nonumber\\
&  -\rho(a,\nabla_{d}^{+}b,c)-\rho(a,b,\nabla_{d}^{+}c). \label{SPS.5a}%
\end{align}
By cycling the letters $a,b,d$ into Eq.(\ref{SPS.5a}), we get
\begin{align}
(\nabla_{a}^{+++-}\rho)(b,d,c)  &  =\nabla_{a}^{+}\rho(b,d,c)-\rho(\nabla
_{a}^{+}b,d,c)\nonumber\\
&  -\rho(b,\nabla_{a}^{+}d,c)-\rho(b,d,\nabla_{a}^{+}c),\label{SPS.5b}\\
(\nabla_{b}^{+++-}\rho)(d,a,c)  &  =\nabla_{b}^{+}\rho(d,a,c)-\rho(\nabla
_{b}^{+}d,a,c)\nonumber\\
&  -\rho(d,\nabla_{b}^{+}a,c)-\rho(d,a,\nabla_{b}^{+}c). \label{SPS.5c}%
\end{align}

Now, by adding Eqs.(\ref{SPS.5a}), (\ref{SPS.5b}), and (\ref{SPS.5c}),
wherever by using Eq.(\ref{TCF.3}) and Eq.(\ref{SPS.2}), we get
\begin{align}
&  (\nabla_{d}^{+++-}\rho)(a,b,c)+(\nabla_{a}^{+++-}\rho)(b,d,c)+(\nabla
_{b}^{+++-}\rho)(d,a,c)\nonumber\\
&  =\nabla_{d}^{+}\rho(a,b,c)+\nabla_{a}^{+}\rho(b,d,c)+\nabla_{b}^{+}%
\rho(d,a,c)\nonumber\\
&  -\rho([a,b],d,c)-\rho([b,d],a,c)-\rho([d,a],b,c)\nonumber\\
&  -\rho(a,b,\nabla_{d}^{+}c)-\rho(b,d,\nabla_{a}^{+}c)-\rho(d,a,\nabla
_{b}^{+}c). \label{SPS.5d}%
\end{align}

But, in accordance with Eq.(\ref{TCF.2a}), we can write
\begin{align}
&  \nabla_{d}^{+}\rho(a,b,c)+\nabla_{a}^{+}\rho(b,d,c)+\nabla_{b}^{+}%
\rho(d,a,c)\nonumber\\
&  =[\nabla_{a}^{+},\nabla_{b}^{+}]\nabla_{d}^{+}c+[\nabla_{b}^{+},\nabla
_{d}^{+}]\nabla_{a}^{+}c+[\nabla_{d}^{+},\nabla_{a}^{+}]\nabla_{b}%
^{+}c\nonumber\\
&  -\nabla_{d}^{+}\nabla_{\lbrack a,b]}^{+}c-\nabla_{a}^{+}\nabla_{\lbrack
b,d]}^{+}c-\nabla_{b}^{+}\nabla_{\lbrack d,a]}^{+}c, \label{SPS.5e}%
\end{align}
and
\begin{align*}
&  -\rho([a,b],d,c)-\rho([b,d],a,c)-\rho([d,a],b,c)\\
&  =-\nabla_{\lbrack a,b]}^{+}\nabla_{d}^{+}c-\nabla_{\lbrack b,d]}^{+}%
\nabla_{a}^{+}c-\nabla_{\lbrack d,a]}^{+}\nabla_{b}^{+}c\\
&  +\nabla_{d}^{+}\nabla_{\lbrack a,b]}^{+}c+\nabla_{a}^{+}\nabla_{\lbrack
b,d]}^{+}c+\nabla_{b}^{+}\nabla_{\lbrack d,a]}^{+}c\\
&  +\nabla_{\lbrack\lbrack a,b],d]}^{+}c+\nabla_{\lbrack\lbrack b,d],a]}%
^{+}c+\nabla_{\lbrack\lbrack d,a],b]}^{+}c.
\end{align*}
i.e., by recalling the Jacobi's identity,
\begin{align}
&  -\rho([a,b],d,c)-\rho([b,d],a,c)-\rho([d,a],b,c)\nonumber\\
&  =-\nabla_{\lbrack a,b]}^{+}\nabla_{d}^{+}c-\nabla_{\lbrack b,d]}^{+}%
\nabla_{a}^{+}c-\nabla_{\lbrack d,a]}^{+}\nabla_{b}^{+}c\nonumber\\
&  +\nabla_{d}^{+}\nabla_{\lbrack a,b]}^{+}c+\nabla_{a}^{+}\nabla_{\lbrack
b,d]}^{+}c+\nabla_{b}^{+}\nabla_{\lbrack d,a]}^{+}c. \label{SPS.5f}%
\end{align}

Now, by adding Eqs.(\ref{SPS.5e}) and (\ref{SPS.5f}), wherever by using
Eq.(\ref{TCF.2a}), we get
\begin{align}
&  \nabla_{d}^{+}\rho(a,b,c)+\nabla_{a}^{+}\rho(b,d,c)+\nabla_{b}^{+}%
\rho(d,a,c)\nonumber\\
&  -\rho([a,b],d,c)-\rho([b,d],a,c)-\rho([d,a],b,c)\nonumber\\
&  =\rho(a,b,\nabla_{d}^{+}c)+\rho(b,d,\nabla_{a}^{+}c)+\rho(d,a,\nabla
_{b}^{+}c). \label{SPS.5g}%
\end{align}

Finally, by putting Eq.(\ref{SPS.5g}) into Eq.(\ref{SPS.5d}), the expected
result immediately follows.$\blacksquare$

\subsection{Cartan Fields}

The \emph{smooth} $(1,2)$-\emph{extensor field} on $U,$ namely $\mathbf{\Theta
},$ which is defined by
\begin{equation}
\mathbf{\Theta}(c)=\frac{1}{2}\partial_{a}\wedge\partial_{b}\tau(a,b)\cdot c
\label{CF.1}%
\end{equation}
will be called the \emph{Cartan torsion field} of $(U,\gamma).$

We should notice that such $\mathbf{\Theta}$ contains all of the geometric
information which is just contained in $\tau.$ Indeed, Eq.(\ref{CF.1}) can be
inverted in such a way that given any $\mathbf{\Theta}$, there is an unique
$\tau$ that verifies Eq.(\ref{CF.1}). We have
\begin{equation}
\tau(a,b)=\partial_{c}(a\wedge b)\cdot\mathbf{\Theta}(c). \label{CF.1a}%
\end{equation}

The \emph{smooth bivector elementary} $2$-\emph{extensor field }on $U$, namely
$\mathbf{\Omega},$ which is defined by
\begin{equation}
\mathbf{\Omega}(c,d)=\frac{1}{2}\partial_{a}\wedge\partial_{b}\rho(a,b,c)\cdot
d \label{CF.2}%
\end{equation}
will be called the \emph{Cartan curvature field }of $(U,\gamma).$

Since Eq.(\ref{CF.2}) can be inverted, by giving $\rho$ in terms of
$\mathbf{\Omega},$ we see that such $\mathbf{\Omega}$ contains the same
geometric information as $\rho.$ The inversion is realized by
\begin{equation}
\rho(a,b,c)=\partial_{d}(a\wedge b)\cdot\mathbf{\Omega}(c,d). \label{CF.2a}%
\end{equation}

\subsection{Cartan's Structure Equations}

Associated to any parallelism structure $(U,\gamma)$ we can introduce two
noticeable operators which map $2$-uples of smooth vector fields to smooth
vector fields. They are:

(a) The mapping $\gamma^{+}:\mathcal{V}(U)\times\mathcal{V}(U)\rightarrow
\mathcal{V}(U)$ defined by
\begin{equation}
\gamma^{+}(b,c)=\partial_{a}(\nabla_{a}^{+}b)\cdot c \label{CSE.1}%
\end{equation}
which will be called the \emph{Cartan connection operator of first kind} of
$(U,\gamma).$

(b) The mapping $\gamma^{-}:\mathcal{V}(U)\times\mathcal{V}(U)\rightarrow
\mathcal{V}(U)$ defined by
\begin{equation}
\gamma^{-}(b,c)=\partial_{a}b\cdot(\nabla_{a}^{-}c) \label{CSE.2}%
\end{equation}
which will be called the \emph{Cartan connection operator of second kind} of
$(U,\gamma).$

We summarize some of the basic properties which are satisfied by the Cartan
operators.\vspace{0.1in}

\textbf{i.} For all $f\in\mathcal{S}(U),$ and $b,b^{\prime},c,c^{\prime}%
\in\mathcal{V}(U)$, we have
\begin{align}
\gamma^{+}(b+b^{\prime},c)  &  =\gamma^{+}(b,c)+\gamma^{+}(b^{\prime
},c),\label{CSE.3a}\\
\gamma^{+}(b,c+c^{\prime})  &  =\gamma^{+}(b,c)+\gamma^{+}(b,c^{\prime
}).\label{CSE.3b}\\
\gamma^{+}(fb,c)  &  =(\partial_{o}f)b\cdot c+f\gamma^{+}(b,c), \label{CSE.3c}%
\\
\gamma^{+}(b,fc)  &  =f\gamma^{+}(b,c). \label{CSE.3d}%
\end{align}

\textbf{ii. }For all $f\in\mathcal{S}(U),$ and $b,b^{\prime},c,c^{\prime}%
\in\mathcal{V}(U)$, we have
\begin{align}
\gamma^{-}(b+b^{\prime},c)  &  =\gamma^{-}(b,c)+\gamma^{-}(b^{\prime
},c),\label{CSE.4a}\\
\gamma^{-}(b,c+c^{\prime})  &  =\gamma^{-}(b,c)+\gamma^{-}(b,c^{\prime
}).\label{CSE.4b}\\
\gamma^{-}(fb,c)  &  =f\gamma^{-}(b,c),\label{CSE.4c}\\
\gamma^{-}(b,fc)  &  =(\partial_{o}f)b\cdot c+f\gamma^{-}(b,c). \label{CSE.4d}%
\end{align}

We have that the Cartan operator of first kind has the linearity property with
respect to the second variable, and the Cartan operator of second kind is
linear with respect to the first variable.

\textbf{iii. }For any $a,b\in\mathcal{V}(U)$,
\begin{equation}
\gamma^{+}(b,c)+\gamma^{-}(b,c)=\partial_{o}(b\cdot c). \label{CSE.5}%
\end{equation}

It is an immediate consequence of Eq.(\ref{CDM.6}).\vspace{0.1in}

\textbf{First Cartan's structure equation}

For any $c\in\mathcal{V}(U)$ it holds
\begin{equation}
\mathbf{\Theta}(c)=\partial_{o}\wedge c+\partial_{s}\wedge\gamma^{-}(s,c).
\label{FCE.1}%
\end{equation}

\textbf{Proof}

By using Eq.(\ref{TCF.1a}) we can write
\begin{align}
\mathbf{\Theta}(c)  &  =\frac{1}{2}\partial_{a}\wedge\partial_{b}(\nabla
_{a}^{+}b-\nabla_{b}^{+}a-[a,b])\cdot c,\nonumber\\
&  =\partial_{a}\wedge\partial_{b}(\nabla_{a}^{+}b-a\cdot\partial_{o}b)\cdot
c. \label{FCE.1a}%
\end{align}

A straightforward calculation yields
\begin{align}
\partial_{a}\wedge\partial_{b}(\nabla_{a}^{+}b)\cdot c  &  =\partial_{a}%
\wedge\partial_{b}\gamma^{+}(b,c)\cdot a\nonumber\\
&  =\partial_{b}\wedge\partial_{a}(\gamma^{-}(b,c)\cdot a-a\cdot\partial
_{o}(b\cdot c)),\nonumber\\
&  =\partial_{b}\wedge\gamma^{-}(b,c)-\partial_{b}\wedge\partial_{o}(b\cdot
c), \label{FCE.1b}%
\end{align}
and
\begin{align}
-\partial_{a}\wedge\partial_{b}(a\cdot\partial_{o}b)\cdot c  &  =\partial
_{a}\wedge\partial_{b}(b\cdot(a\cdot\partial_{o}c)-a\cdot\partial_{o}(b\cdot
c))\nonumber\\
&  =\partial_{a}\wedge(a\cdot\partial_{o}c)+\partial_{b}\wedge\partial
_{a}a\cdot\partial_{o}(b\cdot c),\nonumber\\
&  =\partial_{o}\wedge c+\partial_{b}\wedge\partial_{o}(b\cdot c).
\label{FCE.1c}%
\end{align}

Thus, by putting Eq.(\ref{FCE.1b}) and Eq.(\ref{FCE.1c}) into Eq.(\ref{FCE.1a}%
), we get the expected result.$\blacksquare$\vspace{0.1in}

\textbf{Second Cartan's structure equation}

For any $c,d\in\mathcal{V}(U)$ it holds
\begin{equation}
\mathbf{\Omega}(c,d)=\partial_{o}\wedge\gamma^{+}(c,d)+\gamma^{+}%
(c,\partial_{s})\wedge\gamma^{-}(s,d). \label{SCE.1}%
\end{equation}

\textbf{Proof}

By using Eq.(\ref{TCF.2a}) we have
\begin{align}
\mathbf{\Omega}(c,d)  &  =\frac{1}{2}\partial_{a}\wedge\partial_{b}%
([\nabla_{a}^{+},\nabla_{b}^{+}]c-\nabla_{\lbrack a,b]}^{+}c)\cdot
d,\nonumber\\
&  =\partial_{a}\wedge\partial_{b}(\nabla_{a}^{+}(\nabla_{b}^{+}%
c)-\nabla_{a\cdot\partial_{o}b}^{+}c)\cdot d. \label{SCE.1a}%
\end{align}

But, by taking a pair of reciprocal frame fields $(\{e_{\sigma}\},\{e^{\sigma
}\})$ we can write
\begin{align}
\nabla_{a}^{+}(\nabla_{b}^{+}c)\cdot d  &  =\nabla_{a}^{+}(\gamma
^{+}(c,e^{\sigma})\cdot be_{\sigma})\cdot d\nonumber\\
&  =a\cdot\partial_{o}(\gamma^{+}(c,e^{\sigma})\cdot b)e_{\sigma}\cdot
d+\gamma^{+}(c,e^{\sigma})\cdot b\nabla_{a}^{+}e_{\sigma}\cdot d\nonumber\\
&  =a\cdot\partial_{o}(\gamma^{+}(c,e^{\sigma})\cdot b)e_{\sigma}\cdot
d+\gamma^{+}(c,e^{\sigma})\cdot b\gamma^{+}(e_{\sigma},d)\cdot a,\nonumber\\
&  =a\cdot\partial_{o}\gamma^{+}(c,e^{\sigma})\cdot be_{\sigma}\cdot
d+\gamma^{+}(c,d)\cdot(a\cdot\partial_{o}b)\nonumber\\
&  +\gamma^{+}(c,e^{\sigma})\cdot b\gamma^{+}(e_{\sigma},d)\cdot a.
\label{SCE.1b}%
\end{align}

Now, the first term into Eq.(\ref{SCE.1a}), by using Eq.(\ref{SCE.1b}) and the
well-known identity $\partial_{o}\wedge(fX)=(\partial_{o}f)\wedge
X+f\partial_{o}\wedge X,$ where $f\in\mathcal{S}(U)$ and $X\in\mathcal{M}(U),$
can be written
\begin{align}
\partial_{a}\wedge\partial_{b}\nabla_{a}^{+}(\nabla_{b}^{+}c)\cdot d  &
=\partial_{o}\wedge\gamma^{+}(c,e^{\sigma})(e_{\sigma}\cdot d)+\partial
_{a}\wedge\partial_{b}\gamma^{+}(c,d)\cdot(a\cdot\partial_{o}b)\nonumber\\
&  -\gamma^{+}(c,e^{\sigma})\wedge\gamma^{+}(e_{\sigma},d),\nonumber\\
&  =\partial_{o}\wedge\gamma^{+}(c,e^{\sigma})(e_{\sigma}\cdot d)\nonumber\\
&  -\gamma^{+}(c,e^{\sigma})\wedge(\partial_{o}(e_{\sigma}\cdot d)-\gamma
^{-}(e_{\sigma},d))\nonumber\\
&  +\partial_{a}\wedge\partial_{b}\gamma^{+}(c,d)\cdot(a\cdot\partial
_{o}b),\nonumber\\
&  =\partial_{o}\wedge\gamma^{+}(c,d)\nonumber\\
&  +\gamma^{+}(c,e^{\sigma})\wedge\gamma^{-}(e_{\sigma},d)\nonumber\\
&  +\partial_{a}\wedge\partial_{b}\gamma^{+}(c,d)\cdot(a\cdot\partial_{o}b).
\label{SCE.1c}%
\end{align}

It is also
\begin{equation}
-(\nabla_{a\cdot\partial_{o}b}^{+}c)\cdot d=-\gamma^{+}(c,d)\cdot
(a\cdot\partial_{o}b). \label{SCE.1d}%
\end{equation}

Finally, by putting Eq.(\ref{SCE.1c}) and Eq.(\ref{SCE.1d}) into
Eq.(\ref{SCE.1a}), we get the expected result.$\blacksquare$

\section{Appendix}

Let $U$ be an open subset of $U_{o},$ and let $(U,\phi)$ and $(U,\phi^{\prime
})$ be two local coordinate systems on $U$ compatibles with $(U_{o},\phi
_{o}).$ As we know there must be two pairs of reciprocal frame fields on $U.$
The covariant and contravariant frame fields $\{b_{\alpha}\cdot\partial
x_{o}\}$ and $\{\partial_{o}x^{\alpha}\}$ associated to $(U,\phi),$ and those
ones $\{b_{\mu}\cdot\partial^{\prime}x_{o}\}$ and $\{\partial_{o}%
x^{\mu^{\prime}}\}$ associated to $(U,\phi^{\prime}).$

\textbf{A.1} The $n^{3}$ smooth scalar fields on $U,$ namely $\Gamma
_{\alpha\beta}^{\gamma},$ defined by
\begin{equation}
\Gamma_{\alpha\beta}^{\gamma}=\Gamma^{+}(b_{\alpha}\cdot\partial
x_{o},b_{\beta}\cdot\partial x_{o})\cdot\partial_{o}x^{\gamma} \tag{A1}%
\end{equation}
correspond to the classically so-called \emph{coefficients of connection,} of
course, associated to $(U,\phi).$The coefficients of connection associated to
$(U,\phi^{\prime})$ are given by
\begin{equation}
\Gamma_{\mu^{\prime}\nu^{\prime}}^{\lambda^{\prime}}=\Gamma^{+}(b_{\mu}%
\cdot\partial^{\prime}x_{o},b_{\nu}\cdot\partial^{\prime}x_{o})\cdot
\partial_{o}x^{\lambda^{\prime}}. \tag{A2}%
\end{equation}

We will check next what is the \emph{law of transformation} between them. We
will employ the simplified notations: $b_{\mu}\cdot\partial^{\prime}%
x_{o}=\dfrac{\partial x_{o}^{\sigma}}{\partial x^{\mu^{\prime}}}b_{\sigma}$
and $\partial_{o}x^{\alpha}=b^{\tau}\dfrac{\partial x^{\alpha}}{\partial
x_{o}^{\tau}},$ etc. Then, by recalling the expansion formulas for smooth
vector fields, $v=(v\cdot\partial_{o}x^{\alpha})b_{\alpha}\cdot\partial x_{o}$
and $v=(v\cdot b_{\gamma}\cdot\partial x_{o})\partial_{o}x^{\gamma},$ using
Eqs.(\ref{CO.2a}), (\ref{CO.2b}), (\ref{CO.2c}) and (\ref{CO.2d}), and
recalling the remarkable identity $(b_{\alpha}\cdot\partial x_{o}%
)\cdot\partial_{o}X=b_{\alpha}\cdot\partial X,$ where $X\in\mathcal{M}(U)$
(see \cite{1}), we can write
\begin{align*}
\Gamma_{\mu^{\prime}\nu^{\prime}}^{\lambda^{\prime}}  &  =\Gamma^{+}%
(\frac{\partial x^{\alpha}}{\partial x^{\mu^{\prime}}}b_{\alpha}\cdot\partial
x_{o},\frac{\partial x^{\beta}}{\partial x^{\nu^{\prime}}}b_{\beta}%
\cdot\partial x_{o})\cdot\frac{\partial x^{\lambda^{\prime}}}{\partial
x^{\gamma}}\partial_{o}x^{\gamma}\\
&  =\frac{\partial x^{\alpha}}{\partial x^{\mu^{\prime}}}\frac{\partial
x^{\beta}}{\partial x^{\nu^{\prime}}}\frac{\partial x^{\lambda^{\prime}}%
}{\partial x^{\gamma}}\Gamma_{\alpha\beta}^{\gamma}\\
&  +\frac{\partial x^{\alpha}}{\partial x^{\mu^{\prime}}}(b_{\alpha}%
\cdot\partial x_{o})\cdot\partial_{o}(\frac{\partial x^{\beta}}{\partial
x^{\nu^{\prime}}})\frac{\partial x^{\lambda^{\prime}}}{\partial x^{\gamma}%
}(b_{\beta}\cdot\partial x_{o}\cdot\partial_{o}x^{\gamma})\\
&  =\frac{\partial x^{\alpha}}{\partial x^{\mu^{\prime}}}\frac{\partial
x^{\beta}}{\partial x^{\nu^{\prime}}}\frac{\partial x^{\lambda^{\prime}}%
}{\partial x^{\gamma}}\Gamma_{\alpha\beta}^{\gamma}+\frac{\partial x^{\alpha}%
}{\partial x^{\mu^{\prime}}}\frac{\partial}{\partial x^{\alpha}}%
(\frac{\partial x^{\beta}}{\partial x^{\nu^{\prime}}})\frac{\partial
x^{\lambda^{\prime}}}{\partial x^{\gamma}}\frac{\partial x^{\gamma}}{\partial
x^{\beta}},
\end{align*}
i.e.,
\begin{equation}
\Gamma_{\mu^{\prime}\nu^{\prime}}^{\lambda^{\prime}}=\frac{\partial x^{\alpha
}}{\partial x^{\mu^{\prime}}}\frac{\partial x^{\beta}}{\partial x^{\nu
^{\prime}}}\frac{\partial x^{\lambda^{\prime}}}{\partial x^{\gamma}}%
\Gamma_{\alpha\beta}^{\gamma}+\frac{\partial^{2}x^{\beta}}{\partial
x^{\mu^{\prime}}\partial x^{\nu^{\prime}}}\frac{\partial x^{\lambda^{\prime}}%
}{\partial x^{\beta}}. \tag{A3}%
\end{equation}
It is just the well-known law of transformation for the \emph{classical
coefficients of connection} associated to each of the \emph{coordinate
systems} $\left\langle x^{\mu}\right\rangle $ and $\left\langle x^{\mu
^{\prime}}\right\rangle .\vspace{0.1in}$

\textbf{A.2} Given a smooth vector field on $U,$ say $v,$ the covariant and
contravariant components of $v$ with respect to $(U,\phi)$ and $(U,\phi
^{\prime})$ are respectively given by
\begin{align}
v_{\alpha}  &  =v\cdot(b_{\alpha}\cdot\partial x_{o}),\tag{A4}\\
v^{\alpha}  &  =v\cdot\partial_{o}x^{\alpha},\tag{A5}\\
v_{\alpha}^{\prime}  &  =v\cdot(b_{\alpha}\cdot\partial^{\prime}%
x_{o}),\tag{A6}\\
v^{\alpha^{\prime}}  &  =v\cdot\partial_{o}x^{\alpha^{\prime}}. \tag{A7}%
\end{align}

We will check the relationship between the covariant components of $v$ with
respect to each of the coordinate systems $\left\langle x^{\mu}\right\rangle $
and $\left\langle x^{\mu^{\prime}}\right\rangle .$ By recalling the expansion
formula $w=(w\cdot\partial_{o}x^{\beta})b_{\beta}\cdot\partial x_{o},$ we
have
\begin{align}
v_{\alpha}^{\prime}  &  =v\cdot(b_{\alpha}\cdot\partial^{\prime}x_{o}%
\cdot\partial_{o}x^{\beta})b_{\beta}\cdot\partial x_{o},\nonumber\\
v_{\alpha}^{\prime}  &  =\frac{\partial x^{\beta}}{\partial x^{\alpha^{\prime
}}}v_{\beta}. \tag{A8}%
\end{align}
It is the expected law of transformation for the covariant components of $v$
associated to each of $\left\langle x^{\mu}\right\rangle $ and $\left\langle
x^{\mu^{\prime}}\right\rangle .$

Analogously, by using the expansion formula $w=(w\cdot b_{\beta}\cdot\partial
x_{o})\cdot\partial_{o}x^{\beta},$ we can get the classical law of
transformation for the contravariant components of $v,$ i.e.,
\begin{equation}
v^{\alpha^{\prime}}=\frac{\partial x^{\alpha^{\prime}}}{\partial x^{\beta}%
}v^{\beta}. \tag{A9}%
\end{equation}
\vspace{0.1in}

\textbf{A.3} We will see next which is the meaning of the classical covariant
derivatives of the contravariant and covariant components of a smooth vector
field. By using the expansion formula $v=v^{\alpha}b_{\alpha}\cdot\partial
x_{o},$ and Eqs.(\ref{CDM.4b}) and (\ref{CDM.4c}), and recalling once again
the remarkable identity $(b_{\mu}\cdot\partial x_{o})\cdot\partial_{o}%
X=b_{\mu}\cdot\partial X,$ where $X\in\mathcal{M}(U),$ we can write
\begin{align*}
(\nabla_{b_{\mu}\cdot\partial x_{o}}^{+}v)\cdot\partial_{o}x^{\lambda}  &
=(b_{\mu}\cdot\partial x_{o}\cdot\partial_{o}v^{\alpha})(b_{\alpha}%
\cdot\partial x_{o}\cdot\partial_{o}x^{\lambda})\\
&  +v^{\alpha}(\nabla_{b_{\mu}\cdot\partial x_{o}}^{+}b_{\alpha}\cdot\partial
x_{o})\cdot\partial_{o}x^{\lambda}\\
&  =\frac{\partial v^{\alpha}}{\partial x^{\mu}}\delta_{\alpha}^{\lambda
}+v^{\alpha}\Gamma^{+}(b_{\mu}\cdot\partial x_{o},b_{\alpha}\cdot\partial
x_{o})\cdot\partial_{o}x^{\lambda}\\
&  =\frac{\partial v^{\lambda}}{\partial x^{\mu}}+\Gamma_{\mu\alpha}^{\lambda
}v^{\alpha},
\end{align*}
i.e.,
\begin{equation}
(\nabla_{b_{\mu}\cdot\partial x_{o}}^{+}v)\cdot\partial_{o}x^{\lambda}=\left.
v^{\lambda}\right.  _{;\text{ }\mu}. \tag{A10}%
\end{equation}
By using the expansion formula $v=v_{\alpha}\partial_{o}x^{\alpha},$
Eqs.(\ref{CDM.4b}) and (\ref{CDM.4c}), and Eq.(\ref{CO.3}), we get
\[
(\nabla_{b_{\mu}\cdot\partial x_{o}}^{-}v)\cdot b_{\nu}\cdot\partial
x_{o}=\frac{\partial v_{\nu}}{\partial x^{\mu}}-\Gamma_{\mu\nu}^{\alpha
}v_{\alpha},
\]
i.e.,
\begin{equation}
(\nabla_{b_{\mu}\cdot\partial x_{o}}^{-}v)\cdot b_{\nu}\cdot\partial
x_{o}=\left.  v_{\nu}\right.  _{;\text{ }\mu}. \tag{A11}%
\end{equation}
\vspace{0.1in}

\textbf{A.4} Now, for instance, let us take a smooth $(1,1)$-extensor field on
$U,$ say $t.$ The covariant and contravariant components of $t\ $with respect
to $(U,\phi)$ and $(U,\phi^{\prime})$ are respectively defined to be
\begin{align}
t_{\mu\nu}  &  =t(b_{\mu}\cdot\partial x_{o})\cdot b_{\nu}\cdot\partial
x_{o},\tag{A12}\\
t^{\mu\nu}  &  =t(\partial_{o}x^{\mu})\cdot\partial_{o}x^{\nu},\tag{A13}\\
t_{\mu^{\prime}\nu^{\prime}}  &  =t(b_{\mu}\cdot\partial^{\prime}x_{o})\cdot
b_{\nu}\cdot\partial^{\prime}x_{o},\tag{A14}\\
t^{\mu^{\prime}\nu^{\prime}}  &  =t(\partial_{o}x^{\mu^{\prime}})\cdot
\partial_{o}x^{\nu^{\prime}} \tag{A15}%
\end{align}
It is also possible to introduce two mixed components of $t$ with respect to
$(U,\phi)$ and $(U,\phi^{\prime}).$ They are defined by
\begin{align}
\left.  t_{\mu}\right.  ^{\nu}  &  =t(b_{\mu}\cdot\partial x_{o})\cdot
\partial_{o}x^{\nu},\tag{A16}\\
\left.  t^{\mu}\right.  _{\nu}  &  =t(\partial_{o}x^{\mu})\cdot b_{\nu}%
\cdot\partial x_{o},\tag{A17}\\
\left.  t_{\mu^{\prime}}\right.  ^{\nu^{\prime}}  &  =t(b_{\mu}\cdot
\partial^{\prime}x_{o})\cdot\partial_{o}x^{\nu^{\prime}},\tag{A18}\\
\left.  t^{\mu^{\prime}}\right.  _{\nu^{\prime}}  &  =t(\partial_{o}%
x^{\mu^{\prime}})\cdot b_{\nu}\cdot\partial^{\prime}x_{o}. \tag{A19}%
\end{align}

We will try to check the law of transformation for the covariant component of
$t.$ We can write
\begin{align}
t_{\mu^{\prime}\nu^{\prime}}  &  =t(b_{\mu}\cdot\partial^{\prime}x_{o})\cdot
b_{\nu}\cdot\partial^{\prime}x_{o}\nonumber\\
&  =t(\frac{\partial x^{\alpha}}{\partial x^{\mu^{\prime}}}b_{\alpha}%
\cdot\partial x_{o})\cdot\frac{\partial x^{\beta}}{\partial x^{\nu^{\prime}}%
}b_{\beta}\cdot\partial x_{o}\nonumber\\
&  =\frac{\partial x^{\alpha}}{\partial x^{\mu^{\prime}}}\frac{\partial
x^{\beta}}{\partial x^{\nu^{\prime}}}t(b_{\alpha}\cdot\partial x_{o})\cdot
b_{\beta}\cdot\partial x_{o},\nonumber\\
t_{\mu^{\prime}\nu^{\prime}}  &  =\frac{\partial x^{\alpha}}{\partial
x^{\mu^{\prime}}}\frac{\partial x^{\beta}}{\partial x^{\nu^{\prime}}}%
t_{\alpha\beta}. \tag{A20}%
\end{align}
This is the classical law of transformation for the covariant components of a
smooth $2$-tensor field from $\left\langle x^{\mu}\right\rangle $ to
$\left\langle x^{\mu^{\prime}}\right\rangle .$

By following similar steps we can get the laws of transformation for the
contravariant and mixed components of $t.$ We have
\begin{align}
t^{\mu^{\prime}\nu^{\prime}}  &  =\frac{\partial x^{\mu^{\prime}}}{\partial
x^{\alpha}}\frac{\partial x^{\nu^{\prime}}}{\partial x^{\beta}}t^{\alpha\beta
},\tag{A21}\\
\left.  t_{\mu^{\prime}}\right.  ^{\nu^{\prime}}  &  =\frac{\partial
x^{\alpha}}{\partial x^{\mu^{\prime}}}\frac{\partial x^{\nu^{\prime}}%
}{\partial x^{\beta}}\left.  t_{\alpha}\right.  ^{\beta},\tag{A22}\\
\left.  t^{\mu^{\prime}}\right.  _{\nu^{\prime}}  &  =\frac{\partial
x^{\mu^{\prime}}}{\partial x^{\alpha}}\frac{\partial x^{\beta}}{\partial
x^{\nu^{\prime}}}\left.  t^{\alpha}\right.  _{\beta}. \tag{A23}%
\end{align}
They perfectly agree with the classical laws of transformation for the
respective contravariant and mixed components of a smooth $2$-tensor field
from $\left\langle x^{\mu}\right\rangle $ to $\left\langle x^{\mu^{\prime}%
}\right\rangle $.\vspace{0.1in}

\textbf{A.5} We will get next the relationship among the covariant derivatives
of $t$ and the classical concepts of covariant derivatives of the covariant,
contravariant and mixed components of a smooth $2$-tensor field. For instance,
we have
\begin{align*}
&  (\nabla_{b_{\mu}\cdot\partial x_{o}}^{++}t)(b_{\alpha}\cdot\partial
x_{o})\cdot b_{\beta}\cdot\partial x_{o}\\
&  =b_{\mu}\cdot\partial x_{o}\cdot\partial_{o}(t(b_{\alpha}\cdot\partial
x_{o})\cdot b_{\beta}\cdot\partial x_{o})-t(\nabla_{b_{\mu}\cdot\partial
x_{o}}^{+}b_{\alpha}\cdot\partial x_{o})\cdot b_{\beta}\cdot\partial x_{o}\\
&  -t(b_{\alpha}\cdot\partial x_{o})\cdot\nabla_{b_{\mu}\cdot\partial x_{o}%
}^{+}b_{\beta}\cdot\partial x_{o}\\
&  =\frac{\partial t_{\alpha\beta}}{\partial x^{\mu}}-t(\Gamma^{+}(b_{\mu
}\cdot\partial x_{o},b_{\alpha}\cdot\partial x_{o})\cdot\partial_{o}x^{\sigma
}b_{\sigma}\cdot\partial x_{o})\cdot b_{\beta}\cdot\partial x_{o}\\
&  -t(b_{\alpha}\cdot\partial x_{o})\cdot\Gamma^{+}(b_{\mu}\cdot\partial
x_{o},b_{\beta}\cdot\partial x_{o})\cdot\partial_{o}x^{\tau}b_{\tau}%
\cdot\partial x_{o}\\
&  =\frac{\partial t_{\alpha\beta}}{\partial x^{\mu}}-\Gamma_{\mu\alpha
}^{\sigma}t_{\sigma\beta}-\Gamma_{\mu\beta}^{\tau}t_{\alpha\tau},
\end{align*}
i.e.,
\begin{equation}
(\nabla_{b_{\mu}\cdot\partial x_{o}}^{++}t)(b_{\alpha}\cdot\partial
x_{o})\cdot b_{\beta}\cdot\partial x_{o}=\left.  t_{\alpha\beta}\right.
_{;\text{ }\mu}. \tag{A24}%
\end{equation}
We can also write
\begin{align*}
&  (\nabla_{b_{\mu}\cdot\partial x_{o}}^{+-}t)(b_{\alpha}\cdot\partial
x_{o})\cdot\partial_{o}x^{\beta}\\
&  =b_{\mu}\cdot\partial x_{o}\cdot\partial_{o}(t(b_{\alpha}\cdot\partial
x_{o})\cdot\partial_{o}x^{\beta})-t(\nabla_{b_{\mu}\cdot\partial x_{o}}%
^{+}b_{\alpha}\cdot\partial x_{o})\cdot\partial_{o}x^{\beta}\\
&  -t(b_{\alpha}\cdot\partial x_{o})\cdot\nabla_{b_{\mu}\cdot\partial x_{o}%
}^{-}\partial_{o}x^{\beta}\\
&  =\frac{\partial\left.  t_{\alpha}\right.  ^{\beta}}{\partial x^{\mu}%
}-t(\Gamma^{+}(b_{\mu}\cdot\partial x_{o},b_{\alpha}\cdot\partial x_{o}%
)\cdot\partial_{o}x^{\sigma}b_{\sigma}\cdot\partial x_{o})\cdot\partial
_{o}x^{\beta}\\
&  -t(b_{\alpha}\cdot\partial x_{o})\cdot\Gamma^{-}(b_{\mu}\cdot\partial
x_{o},\partial_{o}x^{\beta})\cdot b_{\tau}\cdot\partial x_{o}\partial
_{o}x^{\tau}.\\
&  =\frac{\partial\left.  t_{\alpha}\right.  ^{\beta}}{\partial x^{\mu}}
-\Gamma_{\mu\alpha}^{\sigma}\left.  t_{\sigma}\right.  ^{\beta}+\Gamma
_{\mu\tau}^{\beta}\left.  t_{\alpha}\right.  ^{\tau},
\end{align*}
i.e.,
\begin{equation}
(\nabla_{b_{\mu}\cdot\partial x_{o}}^{+-}t)(b_{\alpha}\cdot\partial
x_{o})\cdot\partial_{o}x^{\beta}=\left.  \left.  t_{\alpha}\right.  ^{\beta
}\right.  _{;\text{ }\mu}. \tag{A25}%
\end{equation}

\section{Conclusions}

We presented in this paper the notion of a parallelism structure on $M$, i.e.,
a pair $(M,\gamma)$ where $\gamma$ is a \emph{connection }extensor field on
$M$. A theory of $a$-directional covariant derivatives of multivector and
extensor fields is introduced and the main properties satisfied by these
objects are proved. We also give a novel and intrinsic presentation (i.e., one
that does not depend on a chosen orthonormal moving frame) of Cartan theory of
the torsion and curvature fields and of Cartan's structure equations. Examples
as our theory relates to known theories have been worked in detail.

\textbf{Acknowledgments: } V. V. Fern\'{a}ndez and A. M. Moya are very
grateful to Mrs. Rosa I. Fern\'{a}ndez who gave to them material and spiritual
support at the starting time of their research work. This paper could not have
been written without her inestimable help.


\begin{thebibliography}{9}                                                                                                %


\bibitem {1}Moya, A. M., Fern\'{a}ndez, V. V., and Rodrigues, W. A., Jr.
\emph{Multivector and Extensor Fields on Smooth Manifolds}, submitted for
publication \textit{I} (2005).

\bibitem {2}Moya, A. M., Fern\'{a}ndez, V. V., and Rodrigues, W. A., Jr.,
\emph{Extensors in Geometric Algebra}, submitted to publication (2005).

\bibitem {3}Moya, A. M., Fern\'{a}ndez, V. V., and Rodrigues, W. A., Jr.,
\emph{Geometric Algebras}, submitted for publication (2005).
\end{thebibliography}
\end{document}